\documentclass[12pt]{article}
\usepackage{amsmath}
\usepackage{amssymb}
\renewcommand\smallskip{\vskip\smallskipamount}
\renewcommand\medskip{\vskip\medskipamount}
\renewcommand\bigskip{\vskip\bigskipamount}

\begin{document}

\footnotetext{The author is partially supported by NSF Grant
DMS-0707086 and a Sloan Research Fellowship.}

\begin{center}
\begin{large}
\textbf{Nonexistence of Generalized Apparent Horizons in Minkowski
Space}
\end{large}

\bigskip\bigskip
MARCUS A. KHURI
\bigskip\bigskip
\end{center}

\begin{abstract}
We establish a Positive Mass Theorem for initial data sets of the
Einstein equations having generalized trapped surface boundary.
In particular we answer a question posed by R. Wald concerning the
existence of generalized apparent horizons in Minkowski space.
\end{abstract}

\bigskip

\setcounter{equation}{0}

  Let $(M,g,k)$ be an initial data set for the Einstein
equations, that is, $M$ is a Riemannian 3-manifold with metric $g$
and $k$ is a symmetric 2-tensor representing the extrinsic
curvature of a spacelike slice; both are required to satisfy the
constraint equations
\begin{eqnarray*}
16\pi\mu\!\!&=&\!\!R+(\mathrm{Tr}_{g}k)^{2}-|k|^{2},\\
8\pi J_{i}\!\!&=&\!\!\nabla^{j}(k_{ij}-(\mathrm{Tr}_{g}k)g_{ij}),
\end{eqnarray*}
where $R$ is scalar curvature and $\mu$, $J$ are respectively the
energy and momentum densities of the matter fields.  If all
measured energy densities are nonnegative then $\mu\geq|J|$, which
will be referred to as the dominant energy condition.  We assume
that the initial data are asymptotically flat (with one end), so
that at spatial infinity the metric and extrinsic curvature
satisfy the following fall-off conditions
\begin{equation*}
|\partial^{l}(g_{ij}-\delta_{ij})|=O(r^{-l-1}),\text{ }\text{
}\text{ }\text{ }|\partial^{l}k_{ij}|=O(r^{-l-2}),\text{ }\text{
}\text{ }l=0,1,2,\text{ }\text{ }\text{ as }\text{ }\text{
}r\rightarrow\infty.
\end{equation*}
The ADM energy and momentum are then well defined by
\begin{equation*}
E=\lim_{r\rightarrow\infty}\frac{1}{16\pi}\int_{S_{r}}(\partial_{i}g_{ij}-\partial_{j}g_{ii})\nu^{j},\text{
}\text{ }\text{ }\text{
}\overrightarrow{P}_{i}
=\lim_{r\rightarrow\infty}\frac{1}{8\pi}\int_{S_{r}}(k_{ij}-(\mathrm{Tr}_{g}k)g_{ij})\nu^{j},
\end{equation*}
where $S_{r}$ are coordinate spheres in the asymptotic end with
unit outward normal $\nu$.\par
  The strength of the gravitational field in the vicinity of a
2-surface $\Sigma\subset M$ may be measured by the null expansions
\begin{equation*}
\theta_{\pm}:=H_{\Sigma}\pm\mathrm{Tr}_{\Sigma}k,
\end{equation*}
where $H_{\Sigma}$ is the mean curvature with respect to the unit
outward normal (pointing towards spatial infinity). The null
expansions measure the rate of change of area for a shell of light
emitted by the surface in the outward future direction
($\theta_{+}$), and outward past direction ($\theta_{-}$).  Thus
the gravitational field is interpreted as being strong near
$\Sigma$ if $\theta_{+}\leq 0$ or $\theta_{-}\leq 0$, in which
case $\Sigma$ is referred to as a future (past) trapped surface.
Future (past) apparent horizons arise as boundaries of future
(past) trapped regions and satisfy the equation $\theta_{+}=0$
($\theta_{-}=0$).\par
  In an attempt to find the most general conditions under which the
Penrose Inequality is to be valid, Bray and the author [4] have
proposed the notion of a generalized apparent horizon, which we
take to be any surface $\Sigma$ satisfying the equation
\begin{equation*}
H_{\Sigma}=|\mathrm{Tr}_{\Sigma}k|.
\end{equation*}
A very natural question, posed by Wald [10], is to ask whether
such surfaces can exist inside Minkowski space.  Our purpose here
is to show that this is not possible.  The strategy will be to
follow Witten's proof of the Positive Mass Theorem, and show that
if such a surface exists in any initial data set satisfying the
dominant energy condition then the ADM mass is strictly positive,
which of course cannot occur for a slice of Minkowski space.  In
fact this result will be a special case of the Positive Mass
Theorem for spacetimes containing a generalized trapped surface,
that is a surface $\Sigma$ satisfying the inequality
\begin{equation}
H_{\Sigma}\leq |\mathrm{Tr}_{\Sigma}k|.
\end{equation}
It has been shown [5] that the existence of a compact generalized
trapped surface in an asymptotically flat initial data set,
implies the existence of a generalized apparent horizon.  This is
analogous to the relationship between classical trapped surfaces
and apparent horizons [1].  The following theorem exhibits another
analogy between classical and generalized trapped
surfaces.\medskip

\textbf{Theorem.}  \textit{Let $(M,g,k)$ be an asymptotically flat
initial data set for the Einstein equations satisfying the
dominant energy condition $\mu\geq|J|$.  If the boundary $\partial
M$ is nonempty and consists of finitely many compact components
each of which is a generalized trapped surface, then the ADM mass
is strictly positive $E>|\overrightarrow{P}|$.}\medskip


\textit{Proof.}  Let $(\mathcal{M},\gamma)$ be a portion of the
spacetime arising from the initial data $(M,g,k)$, and let
$c:Cl(T\mathcal{M})\rightarrow\mathrm{End}(S)$ be the usual
representation of the Clifford algebra on the bundle of spinors
$S$, so that
\begin{equation*}
c(X)c(Y)+c(Y)c(X)=-2\gamma(X,Y)\mathrm{Id}.
\end{equation*}
We choose a local orthonormal frame $e_{a}$, $a=0,1,2,3$ such that
$e_{0}$ is normal to $M$ and $e_{i}$, $i=1,2,3$ are tangent to
$M$.  Then the ``spacetime spin connection" on $M$ is given by
\begin{equation*}
\nabla_{e_{i}}\psi=(e_{i}(\psi^{I})+\frac{1}{4}\psi^{I}\Gamma_{ij}^{l}c(e^{j})c(e_{l})
+\frac{1}{2}\psi^{I}k_{ij}c(e^{j})c(e_{0}))\phi_{I},
\end{equation*}
where $\psi=\psi^{I}\phi_{I}$ with $\phi_{I}$, $I=1,2,3,4$ being a
choice of spin frame associated with the orhtonormal frame
$e_{i}$, and $\Gamma_{ij}^{l}$ are Christoffel symbols for the
metric $g$.  Note that we are using Dirac spinors $\psi$, which
consist of a pair of $\mathrm{SL}(2,\mathbb{C})$ spinors, one
left-handed and one right-handed.  Consider the following chiral
boundary value problem for the Dirac operator:
\begin{equation}
\mathcal{D}\psi=\sum_{i=1}^{3}c(e_{i})\nabla_{e_{i}}\psi=0\text{
}\text{ }\text{ on }\text{ }\text{ }M,\text{ }\text{
}\psi=\psi_{\infty}+o\left(\frac{1}{|x|^{1-\delta}}\right)\text{
}\text{ }\text{ as }\text{ }\text{ }|x|\rightarrow\infty,
\end{equation}
\begin{equation*}
\epsilon\psi-\psi=0\text{ }\text{ }\text{ on }\text{ }\text{
}\partial M_{+},\text{ }\text{ }\text{ }\text{
}\epsilon\psi+\psi=0\text{ }\text{ }\text{ on }\text{ }\text{
}\partial M_{-},
\end{equation*}
where $\epsilon=c(e_{3})c(e_{0})$ with $e_{3}$ normal to $\partial
M$ and pointing towards spatial infinity, $\psi_{\infty}$ is a
nonzero spinor which is constant in the asymptotic end (the
components $\psi_{\infty}^{I}$, with respect to a fixed frame at
spatial infinity, are constant), and $\partial M_{\pm}$ denotes
the portion of $\partial M$ on which $\theta_{\pm}\leq 0$. Note
that these boundary conditions are the usual ones used to
establish the Positive Mass Theorem with black holes in which
$\partial M$ is assumed to consist of classical future and past
apparent horizons. Below we will show that this boundary value
problem is coercive. Moreover (2) falls into a class of elliptic
boundary value problems treated in [3]. Therefore we conclude that
there exists a unique solution with $\psi-\psi_{\infty}\in
W^{1,2}_{-1}(M)\cap W^{2,2}_{loc}(M)$; here $W^{1,2}_{-1}(M)$ and
$W^{2,2}_{loc}(M)$ represent Sobolev spaces of square integrable
derivatives up to orders one and two respectively, with the
subscript $-1$ indicating an appropriate weight to obtain the
correct fall-off at spatial infinity.\par
  Consider the following Lichnerowicz formula ([3], [7], [9]) for the solution
of (2)
\begin{equation*}
\mathcal{D}^{*}\mathcal{D}\psi=\nabla^{*}\nabla\psi+\mathcal{R}\psi=0,
\end{equation*}
where
\begin{equation*}
\mathcal{R}\psi=4\pi(\mu+J^{i}c(e_{0})c(e_{i}))\psi.
\end{equation*}
Then integrating by parts produces
\begin{eqnarray}
& &\int_{M}\left(|\nabla\psi|^{2}+4\pi(\mu|\psi|^{2}+
J^{i}\langle\psi,c(e_{0})c(e_{i})\psi\rangle)\right)\\
&=&4\pi
P^{a}\langle\psi_{\infty},c(e_{0})c(e_{a})\psi_{\infty}\rangle
-\int_{\partial
M}\langle\psi,c(e_{3})\sum_{i=1}^{2}c(e_{i})\nabla_{e_{i}}\psi\rangle,\nonumber
\end{eqnarray}
where $P^{a}$ is the ADM 4-momentum.  In order to facilitate
calculation of the boundary term, we define the boundary covariant
derivative by
\begin{equation*}
\overline{\nabla}_{e_{i}}\psi=e_{i}(\psi)+\frac{1}{4}\sum_{j,l=1}^{2}\Gamma_{ij}^{l}
c(e_{j})c(e_{l})\psi+\frac{1}{2}k_{i3}c(e_{3})c(e_{0})\psi,\text{
}\text{ }\text{ }\text{ }i=1,2,
\end{equation*}
and we define the boundary Dirac operator by
\begin{equation*}
\mathcal{D}_{\partial
M}\psi=c(e_{3})\sum_{i=1}^{2}c(e_{i})\overline{\nabla}_{e_{i}}\psi.
\end{equation*}
The boundary term may now be calculated by using properties of
Clifford multiplication, symmetries of the Christoffel symbols,
and the special boundary conditions of (2), as follows
\begin{eqnarray}
&
&c(e_{3})\sum_{i=1}^{2}c(e_{i})\nabla_{e_{i}}\psi\\
&=&\mathcal{D}_{\partial
M}\psi+\frac{1}{2}\sum_{i=1}^{2}\sum_{j=1}^{3}\Gamma_{ij}^{3}c(e_{3})c(e_{i})
c(e_{j})c(e_{3})\psi
+\frac{1}{2}\sum_{i,j=1}^{2}k_{ij}c(e_{3})c(e_{i})c(e_{j})c(e_{0})\psi\nonumber\\
&=&\mathcal{D}_{\partial M}\psi-\frac{1}{2}H_{\partial M}\psi
-\frac{1}{2}(\mathrm{Tr}_{\partial M}k)c(e_{3})c(e_{0})\psi\nonumber\\
&=&\mathcal{D}_{\partial M}\psi-\frac{1}{2}\theta_{\pm}\psi\text{
}\text{ }\text{ on }\text{ }\text{ }\partial M_{\pm}.\nonumber
\end{eqnarray}
Moreover a similar calculation shows that $\mathcal{D}_{\partial
M}\epsilon=-\epsilon\mathcal{D}_{\partial M}$, and therefore since
$\epsilon$ is self-adjoint with respect to the (positive definite)
inner product $\langle\cdot,\cdot\rangle$ on $S$ we have
\begin{eqnarray}
\langle\psi,\mathcal{D}_{\partial
M}\psi\rangle&=&\pm\langle\psi,\mathcal{D}_{\partial
M}\epsilon\psi\rangle\\
&=&\mp\langle\psi,\epsilon\mathcal{D}_{\partial
M}\psi\rangle\nonumber\\
&=&\mp\langle\epsilon\psi,\mathcal{D}_{\partial
M}\psi\rangle=-\langle\psi,\mathcal{D}_{\partial
M}\psi\rangle\nonumber.
\end{eqnarray}
Then by combining (3), (4), (5), choosing $\psi_{\infty}$ so that
\begin{equation*}
P^{a}\langle\psi_{\infty},c(e_{0})c(e_{a})\psi_{\infty}\rangle =
E-|\overrightarrow{P}|,
\end{equation*}
and applying the dominant energy condition, it follows that
\begin{equation}
\int_{M}|\nabla\psi|^{2}\leq 4\pi(E-|\overrightarrow{P}|).
\end{equation}
Note that the same arguments used to obtain this inequality also
yield the coercivity of boundary value problem (2), which is
needed for establishing the existence and regularity of
solutions.\par
  We now proceed by contradiction and assume that $E\leq |\overrightarrow{P}|$.
Then (6) shows that $\psi$ is covariantly constant.  First
consider the case in which at least one boundary component
$\Sigma$ is a true generalized trapped surface, that is $\Sigma$
satisfies (1) and $\mathrm{Tr}_{\Sigma}k$ changes sign along
$\Sigma$. Then according to the boundary conditions imposed on
$\psi$, and the fact that $\psi$ is continuous up to the boundary
(with the help of a Sobolev embedding), there is a point
$p\in\Sigma$ at which $\psi(p)=0$.  Now parallel transport $\psi$
along any curve emanating from $p$ to find that $\psi=0$ along
this curve (since $\psi$ restricted to the curve is itself the
solution of parallel transport).  But this implies that
$\psi\equiv 0$ on $M$, which is impossible as $\psi_{\infty}\neq
0$.\par
  In the remaining case to consider all boundary components are either
pure future or past trapped surfaces.
Then according to Andersson and Metzger [1] there exists a smooth
compact outermost apparent horizon, each component of which either
has spherical topology or is a flat torus [6].  First assume that
at least one component $\Sigma$ of the outermost apparent horizon
has spherical topology, and we further assume that it is a future
apparent horizon, that is $\theta_{+}=0$ (similar arguments will
hold for a past apparent horizon).  By writing the full Dirac
operator in terms of the induced operator on the boundary with the
help of calculation (4), and using the fact that $\psi$ is
covariantly constant, on $\Sigma$ we have
\begin{equation*}
0=c(e_{3})\mathcal{D}\psi=-\nabla_{e_{3}}\psi+\mathcal{D}_{\Sigma}\psi
-\frac{1}{2}\theta_{+}\psi=\mathcal{D}_{\Sigma}\psi.
\end{equation*}


\noindent Since $\psi$ cannot vanish on $\Sigma$ (according to
arguments above), this says that $\Sigma$ admits a nontrivial
harmonic spinor.  However since $\Sigma$ is topologically a
2-sphere, this is impossible according to the Hijazi-B\"{a}r
inequality ([2], [8]), which states that all eigenvalues of the
Dirac operator on a 2-sphere must satisfy
\begin{equation*}
|\lambda(\mathcal{D}_{\Sigma})|\geq\sqrt{\frac{4\pi}{\mathrm{Area}(\Sigma)}}>0.
\end{equation*}
Application of the Hijazi-B\"{a}r inequality was first suggested
by Bartnik and Chru\'{s}ciel in [3].\par
  If all components of the outermost apparent horizon are flat tori
then we proceed as follows.  Assume that $\partial M$ coincides
with the outermost apparent horizon. Then we may choose a spin
structure on $M$ for which the induced spin structure on one of
the boundary components $\Sigma$ is not the ``trivial" one (note
that the 2-torus admits four distinct spin structures).  By
arguing as above we find that $\Sigma$ admits a nontrivial
harmonic spinor.  However this is impossible, since only the
trivial spin structure on a flat torus admits nontrivial harmonic
spinors [2].  With this contradiction we conclude that the ADM
mass must be strictly positive.  Q.E.D.\medskip

\begin{center}
\textbf{References}
\end{center}
\medskip

\noindent 1.\hspace{.03in} L. Andersson, and J. Metzger,
\textit{The area of horizons and the trapped region},
\par\hspace{-.05in} preprint, arXiv:0708.4252, 2007.\medskip

\noindent 2.\hspace{.03in} C. B\"{a}r, \textit{Harmonic spinors
and topology}, New Developments in Differential Geom-
\par\hspace{-.05in} etry, Budapest 1996, 53-66.\medskip

\noindent 3.\hspace{.03in} R. Bartnik, and P. Chru\'{s}ciel,
\textit{Boundary value problems for Dirac-type equations},
\par\hspace{-.05in} J. Reine Angew. Math., $\mathbf{579}$ (2005), 13-73.\medskip

\noindent 4.\hspace{.03in} H. Bray, and M. Khuri, \textit{PDE's
which imply the Penrose conjecture}, in preparation,
\par\hspace{-.05in} 2008.\medskip

\noindent 5.\hspace{.03in} M. Eichmair, \textit{Existence,
regularity, and properties of generalized apparent horizons},
\par\hspace{-.05in} preprint, arXiv:0805.4454, 2008.\medskip

\noindent 6.\hspace{.03in} G. Galloway, and R. Schoen, \textit{A
generalization of Hawking's black hole topology}
\par\hspace{-.05in} \textit{theorem to higher dimensions},
Comm. Math. Phys., $\mathbf{266}$ (2006), no. 2, 571-576.\medskip

\noindent 7.\hspace{.03in} M. Herzlich, \textit{The positive mass
theorem for black holes revisited},  J. Geom. Phys.,
\par\hspace{-.05in} $\mathbf{26}$ (1998), no. 1-2, 97-111.\medskip

\noindent 8.\hspace{.05in}  O. Hijazi, \textit{Premi\`{e}re valeur
propre de l'op\'{e}rateur de Dirac et nombre de Yamabe},
\par\hspace{-.05in} C.R. Acad. Sci. Paris, $\mathbf{313}$ (1991),
865-868.\medskip

\noindent 9.\hspace{.06in} T. Parker, and C. Taubes, \textit{On
Witten's proof of the positive energy theorem},
\par\hspace{-.03in} Commun. Math. Phys., $\mathbf{84}$ (1982), no. 2, 223-238.\medskip

\noindent 10.\hspace{-.03in}  \textit{Workshop on Mathematical
Aspects of General Relativity}, Niels Bohr Interna-
\par\hspace{-.01in} tional Academy (Copenhagen, Denmark), April
7-17, 2008.

\bigskip\bigskip\footnotesize

\noindent\textsc{Department of Mathematics, Stony Brook
University, Stony Brook, NY 11794}\par

\noindent\textit{E-mail address}: \verb"khuri@math.sunysb.edu"

\end{document}